\documentclass[11pt]{article}
\usepackage{fullpage,amsmath,amsthm,amsfonts,verbatim,amssymb,epsfig,txfonts}
\usepackage{hyperref}

\DeclareMathOperator{\E}{\mathbb{E}}

\newcommand{\R}{\mathbb{R}}

\newcommand{\Z}{\mathbb{Z}}

\newcommand{\Le}{\Biggl}
\newcommand{\Ri}{\Biggr}

\newcommand{\e}{\varepsilon}

\theoremstyle{plain}
\newtheorem{theorem}{Theorem}

\newtheorem{corollary}[theorem]{Corollary}

\theoremstyle{definition}

\begin{document}

\title{An application of metric cotype to quasisymmetric embeddings}

\date{}

\author{ Assaf Naor}

\maketitle

\begin{abstract}
We apply the notion of metric cotype to show that $L_p$ admits a
quasisymmetric embedding into $L_q$ if and only if  $p\le q$ or
$q\le p\le 2$.
\end{abstract}

This note is a companion to~\cite{MN05-journal}. After the final
version of~\cite{MN05-journal} was sent to the journal for
publication I learned from Juha Heinonen and Leonid Kovalev of a
long-standing open problem in the theory of quasisymmetric
embeddings, and it turns out that this problem can be resolved using
the methods of~\cite{MN05-journal}. The argument is explained below.
I thank Juha Heinonen and Leonid Kovalev for bringing this problem
to my attention.

Let $(X,d_X)$ and $(Y,d_Y)$ be metric spaces. An embedding $f:X\to
Y$ is said to be a quasisymmetric embedding with modulus $\eta:
(0,\infty)\to (0,\infty)$ if $\eta$ is increasing, $\lim_{t\to
0}\eta(t)=0$, and for every distinct $x,y,z\in X$ we have
$$
\frac{d_Y(f(x),f(y))}{d_Y(f(x),f(z))}\le \eta
\left(\frac{d_X(x,y)}{d_X(x,z)}\right).
$$
We refer to~\cite{Heinonen01} and the references therein for a
discussion of this notion.

It was not known whether every two separable Banach spaces are
quasisymetrically equivalent. This is asked in~\cite{Vai99} (see
problem 8.3.1 there). We will show here that the answer to this
question is negative. Moreover, it turns out that under mild
assumptions the cotype of a Banach space is preserved under
quasisymmetric embeddings. Thus, in particular, our results imply
that $L_p$ does not embed quasisymetrically into  $L_q$ if $p>2$ and
$q<p$. The question of determining when $L_p$ is quasisymetrically
equivalent to $L_q$ was asked in~\cite{Vai99} (see problem 8.3.3
there). We also deduce, for example, that the separable space $c_0$
does not embed quasisymetrically into any Banach space which has an
equivalent uniformly convex norm.

We recall some definitions. A Banach space $X$ is said to have
(Rademacher) type $p> 0$ if there exists a constant $T<\infty$ such
that for every $n$ and every $x_1,\ldots,x_n\in X$,
\begin{eqnarray*}\label{eq:def type}
\E_\e\Le\|\sum_{j=1}^n \e_j x_j\Ri\|_X^p\le T^p\sum_{j=1}^n
\|x_j\|_X^p.
\end{eqnarray*}
where the expectation $\E_\e$ is with respect to a uniform choice of
signs $\e=(\e_1,\ldots,\e_n)\in \{-1,1\}^n$. $X$ is said to have
(Rademacher) cotype $q>0$ if there exists a constant $C<\infty$ such
that for every $n$ and every $x_1,\ldots,x_n\in X$,
\begin{eqnarray*}\label{eq:def Rademacher cotype}
\E_\e\Le\|\sum_{j=1}^n \e_j x_j\Ri\|_X^q\ge
\frac{1}{C^q}\sum_{j=1}^n \|x_j\|_X^q.
\end{eqnarray*}
We also write
$$
p_X=\sup\{p\ge 1:\ X\ \text{has\  type}\  p\}\quad \mathrm{and}\quad
q_X=\inf\{q\ge 2:\ X\ \text{has\  cotype}\  q\}.
$$
$X$ is said to have non-trivial type if $p_X>1$, and $X$ is said to
have non-trivial cotype if $q_X<\infty$. For example, $L_p$ has type
$\min\{p,2\}$ and cotype $\max\{p,2\}$ (see for example~\cite{MS86}.

\begin{theorem}\label{thm:quasi} Let $X$ be a
Banach space with non-trivial type. Assume that $Y$ is a Banach
space which embeds quasisymmetrically into $X$. Then $q_Y\le q_X$.
\end{theorem}

\begin{proof} Let $f:Y\to X$ be a quasisymmetric embedding with modulus $\eta$. Assume for the sake of contradiction that $X$ has cotype $q$ and that
$p\coloneqq q_Y>q$. By the Maurey-Pisier theorem~\cite{MP76} for
every $n\in \mathbb N$ there is a linear operator $T:\ell_p^n\to Y$
such that for all $x\in \ell_p^n$ we have $\|x\|_p\le \|T(x)\|_Y\le
2\|x\|_p$. For every integer $m\in \mathbb N$ consider the mapping
$g: \Z_m^n\to X$ given by
$$
g(x_1,\ldots,x_n)=f\circ T\Le(e^{\frac{2\pi i
x_1}{m}},\ldots,e^{\frac{2\pi i x_n}{m}}\Ri).
$$
By Theorem 4.1 in~\cite{MN05-journal} there exist constants $A,B>0$
which depend only on the type and cotype constants of $X$ such that
for every integer $m\ge An^{1/q}$ which is divisible by $4$ and
every $h:\Z_m^n\to X$ we have
\begin{eqnarray}\label{eq:def cotype}
\sum_{j=1}^n\E_x\left[\left\|h\left(x+\frac{m}{2}e_j\right)-h(x)\right\|_X^q\right]\le
B^q m^q\E_{\e,x}\left[\|h(x+\e)-h(x)|_X^q\right],
\end{eqnarray}
where the expectations above are taken with respect to uniformly
chosen $x\in \Z_m^n$ and $\e\in\{-1,0,1\}^n$ (here, and in what
follows we denote by $\{e_j\}_{j=1}^n$ the standard basis of
$\R^n$).

From now on we fix $m$ to be be the smallest integer which is
divisible by $4$ and $m\ge An^{1/q}$. Thus $m\le 8An^{1/q}$. For
every $x\in \Z_m^n$, $j\in \{1,\ldots,n\}$ and $\e\in \{-1,0,1\}^n$
we have
$$
\frac{\|g(x+\e)-g(x)\|_X}{\left\|g\left(x+\frac{m}{2}e_j\right)-g(x)\right\|_X}\le
\eta\left(\frac{\left\|T\left(\sum_{k=1}^n \left( e^{\frac{\pi i
\e_k}{m}}-1\right)e_j\right)\right\|_Y}{\|T(2e_j)\|_Y}\right)\le
\eta\left(\frac{\pi n^{1/p}}{m}\right)\le  \eta\left(\frac{\pi}{A}
n^{\frac{1}{p}-\frac{1}{q}}\right).
$$
Thus, using~\eqref{eq:def cotype} for $g=h$ we see that
\begin{multline*}
n\E_{\e,x}\|g(x+\e)-g(x)\|_X^q\le \eta\left(\frac{\pi}{A}
n^{\frac{1}{p}-\frac{1}{q}}\right)^q\sum_{j=1}^n\E_x\left\|g\left(x+\frac{m}{2}e_j\right)-g(x)\right\|_X^q\\\le
\eta\left(\frac{\pi}{A} n^{\frac{1}{p}-\frac{1}{q}}\right)^q(8AB)^qn
\E_{\e,x}\|g(x+\e)-g(x)\|_X^q.
\end{multline*}
Canceling the term $\E_{\e,x}\|g(x+\e)-g(x)\|_X^q$ we deduce that
$$
\eta\left(\frac{\pi}{A} n^{\frac{1}{p}-\frac{1}{q}}\right)\ge
\frac{1}{8AB}.
$$
Since $p>q$ this contradicts the fact that $\lim_{t\to 0}
\eta(t)=0$.
\end{proof}

Using the same argument as in~\cite{MN05-journal} (and noting that
the snowflake embedding from~\cite{MN04} is a quasisymmetric
embedding), we obtain the following complete answer to the question
when $L_p$ embeds quasisymmetrically into $L_q$.

\begin{corollary}\label{thm:quasiL_p}
For $p,q> 0$, $L_p$ embeds quasisymmetrically into $L_q$ if and
only if $p\le q$ or $q\le p\le 2$.
\end{corollary}

%\bibliographystyle{abbrv}

%\bibliography{cotype}

\def\cprime{$'$}

%{\small
%\bibliographystyle{abbrv}
%\bibliography{cotype}
%}

\end{document}